\documentclass[%
leqno,%
12pt%
]{article}

\usepackage[utf8]{inputenc}
\usepackage[T1]{fontenc}
\usepackage[english]{babel}

\usepackage[intlimits]{amsmath}
\usepackage{amssymb, amsfonts}
\usepackage[amsmath, thmmarks]{ntheorem}

\usepackage{mathrsfs}
\usepackage{bbm}
\usepackage{graphicx}
\usepackage{float}

\usepackage{booktabs}

\newcounter{theoremcounter}

\theoremstyle{plain}
\theoremseparator{.}

\theoremstyle{plain}

\newtheorem{theorem}[theoremcounter]{Theorem}

\newtheorem{conjecture}[theoremcounter]{Conjecture}

\theorembodyfont{\normalfont}
\theoremsymbol{\ensuremath{\Box}}

\newcommand{\NN}{\ensuremath{\mathbb{N}}}
\newcommand{\ZZ}{\ensuremath{\mathbb{Z}}}
\newcommand{\QQ}{\ensuremath{\mathbb{Q}}}
\newcommand{\RR}{\ensuremath{\mathbb{R}}}
\newcommand{\CC}{\ensuremath{\mathbb{C}}}

\newcommand{\tr}{\ensuremath{\mathrm{tr}}}

\newcommand{\SL}[2]{\ensuremath{\mathrm{SL}_{#1}(#2)}}

\newcommand{\isdiv}{\ensuremath{\mid}}

\title{Efficiently generated spaces of classical Siegel modular forms and the Böcherer conjecture}
\author{Martin Raum\\\small{MRaum@mpim-bonn.mpg.de}\vspace{1.5em}\\ MPI für Mathematik\thanks{major parts this the work have been done at Lehrstuhl A f\"ur Mathematik, \mbox{RWTH Aachen University}, 52056 Aachen, Germany}\\Vivatsgasse 7\\53111 Bonn, Germany}
\date{\today}

\begin{document}
\maketitle

\begin{abstract}
We state and verify up to weight $172$ a conjecture on the existence of a certain generating set for spaces of classical Siegel modular forms.
This conjecture is particularly useful for calculations involving Fourier expansions.
Using this generating set we verify the Böcherer conjecture for non-rational eigenforms.
As one further application we verify another conjectures for weights up to $150$ and investigate an analogue of the Victor-Miller basis.
Additionally, we describe some arithmetic properties of the basis we found.
\end{abstract}

\section*{Introduction}

Nowadays, we know many basic facts about Siegel modular forms of degree $2$.
In particular, the classical Siegel modular forms were investigated thoroughly.
However, we still lack knowlegde about analogues of well known properties of elliptic modular forms.
A first effort to overcome this situation by a computational approach has been done in \cite{skoruppasiegelcomputations}.
Moreover, there are recent calculations of Hecke eigenvalues of vector valued Siegel modular forms by Bergström, Carel and van de Geer \cite{cohomologiacalcomputations} using cohomological methods.

Both computations have in common that they restrict to low weights. Moreover, they focus on Fourier expansions. But neither of them is able to provide access to the associated L-series up to sufficient precision. This shows that we need an efficient way to generate the Fourier expansions of the spaces for weight $k$ forms, if $k$ is big.

We provide evidence for a conjecture on generating sets for spaces of classical Siegel modular forms.
Namely, the products of at most two elements of the Maass spezialschar form a generating set of the space of all Siegel modular forms associated to the full modular group.
We will use Sage\cite{sage} for our computations. Additionally, we will use Magma\cite{magma} for fast linear algebra.

In section \ref{sec:generatingsets} we will verify up to weight $172$ that these product generate the spaces of Siegel modular forms.
This makes it possible to calculate Fourier expansions for very high weight forms up to high precisions.
These possiblities can be used to compute the attached spinor L-series of Hecke eigenforms up to a precision involving primes up to $1500$.
\vspace{1.5em}

We will do this in section \ref{sec:boechererverification} to investigate the following conjecture.
In \cite{boechererconjecture} an analog of the Gross-Zagier theorem for the spinor L-series attached to Siegel modular forms of degree $2$ has been conjectured.
\begin{conjecture}[Böcherer]
\label{boechererconjecture_explicit}
Let $f(z) = \sum_t a_f(t) \exp(2 \pi i \cdot \tr(t z))$ be a Siegel Hecke eigenform of degree $2$ associated to the full modular group. Denote the attached spinor L-series twisted by the Kronecker character $\left( \frac{D}{\cdot}\right),\, D < 0$ by $Z_{f, D}$. Then
\begin{align*}
\label{boechererconjectureformula}
  Z_{f,D} (k - 1)
&=
  c_f B_f (D) \qquad \text{with}
\\
  B_f (D)
&=
  \left( \sum_{[t]\, \text{s.t.}\, -4\cdot\det t = D} \hspace{-0.2em}
   a_f(t) / |\mathrm{Aut}(t)| \right)^2
\end{align*}
for some $c_f \in \CC$ depending only on $f$. Here $[t]$ runs over a set of representatives with respect to the $\SL{2}{\ZZ}$ action and $\mathrm{Aut} (t) = \{a \in \SL{2}{\ZZ} \,:\,a^\tr t a = t\}$.
\end{conjecture}

An affirmative treatment of this conjecture could also provide a tool to investigate statistics of central values more efficiently.

The conjecture has been proven for Maass lifts (cf. \cite{boechererconjecture}) by using the splitting of the L-series and Waldspurger's theorem on special values of L-series attached to elliptic modular forms.
For non Maass lifts it seems inaccessible at the moment. Indeed, even the functional equation of the twisted L-series, which appears in the conjecture, has be proven completely only recently \cite{spinorfunctionalequation}.

The conjecture has been verified in \cite{kohnenkuss_boechererconjecture} for even weights $k$ from $20$ to $26$ and discriminants $D \in \{-3,-4,-7,-8\}$.
Besides the small number of cases considered, there is another outstanding problem with this verification.
Up to weight $26$ all Hecke eigenforms are rational. But the space of cuspidal non Maass lifts has dimension $2$ for weight $24$ and $26$. Hence, the rational Hecke action on this space is not irreducible.
This phenomena is exceptional as suggested by computations presented in section \ref{sec:applications}.
Namely, the rational action of the Hecke algebra on spaces of weight $k$ Siegel modular forms is irreducible if $28 \le k \le 150$.
The reason for this curiosity is not known yet, but checking the Böcherer conjecture outside the range of this phenomena can provide further support for it.

This is done in section \ref{sec:boechererverification}. Namely, we optain intervals $I_{f,D}$ such that $c_f$ (cf. conjecture \ref{boechererconjectureformula}) is contained in $I_{f,D}$, if the conjecture holds. The intersection $\cap_D I_{f,D}$ is non empty for weight $k$ cuspidal Hecke eigenforms which are not contained in the Maass spezialschar, $20 \le k < 40$.
\vspace{1.5em}

In section \ref{sec:applications} we present two further applications . Firstly, we state a conjecture on the maximal discriminant of a pivot set of a basis' Fourier expansion. In a recent preprint \cite{pooryuenparamoduluar} this has already been considered. We are able to provide evidence of a better asympthotic behaviour.
Secondly, we consider the rational Hecke action on the spaces of weight $k$ forms as mentioned above.

Although the construction of the generating set was purely motivated by computational needs, it turned out to have intesting arithmetic properties.
We will present these at the end of section \ref{sec:applications}.

Throughout the whole paper we only consider spaces of rational modular forms. In particular we are not concerned with integral modular forms.
A basis of the module of integral Siegel modular forms can easily be deduced and might be of some use, too.

\subsubsection*{Acknowledgements}
\textit{The author thanks the Lehrstuhl B für Mathematik, RWTH Aachen University, 52056 Aachen, Germany for the computational resources they have provided.
}

\section{Generating spaces by products of elements in the Maass spezialschar}
\label{sec:generatingsets}

The well known Igusa generators for the ring of classical Siegel modular forms associated to the full modular group of degree $2$  (cf. \cite{igusagenerators}) are elements of the Maass spezialschar.
The multiplication of multi dimensional Fourier expansions is very expensive. Thus, an important question raises. Namely, how many multiplications will we need to calculate a basis element of a given space of modular forms of weight $k$?
Using Igusa generators the answer is clearly $O(\log(k))$.
For example, calculating the Fourier expansion of a generic weight $30$ form, we have to perform more than $15$ multiplications and the longest product involves $4$ multiplications.

Since it is rather cheap to calculate Fourier expansions of elements in the Maass spezialschar we do not need to restrict to the aforementioned generators of the ring, but we can consider all elements of this space.
Indeed, up to weight $18$ the Maass spezialschar equals the spaces of modular forms. For every higher weight it is a proper subspace.
Surprisingly, using this enlarged set of generators we are able to restrict to products with at most 2 factors.
Using an according basis to calculate the fourier expansion of a generic weight $30$ form, we need $6$ multiplications and the longest product involves two factors.

To state a precise conjecture denote the Jacobian modular group of degree $1$ by $\Gamma_J$ and the full Siegel modular group of degree $g$ by $\Gamma_g$. Moreover let $[\Gamma_g, \det^k]$ and $[\Gamma_g, \det^k]_0$ be the spaces of modular forms and cusp forms of weight $k$, respectively. In analogy denote the space of Jacobi forms of index $m$ and weight $k$ by $[\Gamma_J, m, k]$.
Let $V f$ denote Gritsenko's arithmetic lift of a Jacobi form $f$ (cf. \cite{gritsenkoparamodularforms}).
\begin{conjecture}\label{generatingspaces}
We have
\begin{align*}
  \bigoplus_{k \ge 0,\, 2 \isdiv k} [\Gamma_2, {\det}^k]
&= 
  \left( \bigoplus_{k \ge 0,\, 2 \isdiv k} V [\Gamma_J, 1, k] \right)^2
.
\end{align*}
\end{conjecture}

By explicit computations we are able to proof the following theorem.
\begin{theorem}
Conjecture \ref{generatingspaces} holds up to weight $172$.
\end{theorem}

In table \ref{spanningmaassproductstable} we present some of the generators that we calculated. We fix some notation. The Maass spezialschar is the image of $[\Gamma_J,1,k] \cong [\Gamma_1, \det^k] \times [\Gamma_1, \det ^{k+1}]_0$. Every space $[\Gamma_1, \det^k]$  admits a Victor-Miller basis, which is canonically enumerated. We denote this basis by $f_{k,\,i}$ for $i \in \{0,\ldots,\dim [\Gamma_1, \det ^k] - 1\}$. The lift of $(f_{k,\,i},0)$ is denoted by $k |_i$ and the lift of $(0, f_{k+2,i})$ is denoted by $k|^i$. In the table we have listed only elements which do not belong to the Maass spezialschar. Together with $V [\Gamma_J, 1, k]$ they span $[\Gamma_2, \det ^k]$.


\section{The Böcherer conjecture}
\label{sec:boechererverification}

To verify the Böcherer conjecture for non-rational eigenforms a suitable generating set for the spaces of Siegel modular forms was needed. This motivated the investigations in section \ref{sec:generatingsets}.

In this section we will approximate the central values of the twisted spinor \mbox{L-function} with a bound on the approximation error.
This yields very small intervalls for the constant $c_f$.
Since they have non empty intersection, the calculations in this section support conjecture \ref{boechererconjecture_explicit}.

Consider some cuspidal eigenform $f \in [\Gamma_2, \det ^k]$ and a primitive Dirichlet character $\chi_D = \left( \frac{D}{\cdot} \right)$.
We revise the approximation of the twisted spinor \mbox{L-function} $Z_{f, D}$ at some point $s \in \CC$ based on some series representation found by Kohnen (cf. \cite{kohnenkuss_boechererconjecture}).
We denote the n-th coefficient of the spinor \mbox{$\zeta$-function} by $\lambda_{f} (n)$. Then
\begin{align*}
  Z_{f,D} (s)
& =
  \sum_{n = 1} ^\infty g_D (n) \lambda_f (n) \chi_D(n) \qquad\text{with}
\\
  g_D (n)
& =
  2 (2 \pi / |D|)^{2 - k + 2 s}
  n^{-k/2} / (\Gamma(s) \Gamma(s - k + 2))
\\&\qquad
  \cdot
  \int_n ^\infty K_{k-2} (4 \pi \sqrt y / |D|)
                 \left((y/n)^{3k/2 - 2 - s} + (y/n)^{s - k/2} \right) dy
.
\end{align*}

This series can be approximated by truncating it suitably.
\begin{align*}
  \widetilde{Z}_{f,D} (s, P, N)
=
  \sum_{\substack{1\le n \le N \,\text{s.t.}\\ \forall p\,\text{prime} \isdiv n \,:\, p < P}}
  \hspace{-1.5em}
  g_D (n) \lambda_f (n) \chi_D(n)
.
\end{align*}

If $s = k - 1$, a bound for the approximation error has already been given by Kohnen and Kuß.
For general $s = \sigma + i t,\,\sigma,t \in \RR$ an easy calculation yields the following formula, which is numerically stable with respect to $\sigma$.
\begin{align*}
&\qquad
  \left| Z_{f,D}(s) - \widetilde{Z}_{f,D}(s,P,N) \right| 
\le
  \eta_f(N,P,D,s)
\\&
:=
  \sum_{\substack{N > p > P\\1 \le \nu \le N/p}}
  4 p^ {k - 3/2} |\lambda_f (\nu)| \widetilde{g}_D (\nu p)
\\&\quad
  + 2 \left(2 \pi / |D|\right)^{2 - k + 2 s} / (\Gamma(s) \Gamma(s - k + 2))
\\&\quad
  \cdot \Bigg(
  \left(\mathbbm{1}_{\sigma > k - 5/2} + \mathbbm{1}_{\sigma < k + 1/2}\right)
  \int_{N+1} ^\infty K_{k-2} (4 \pi \sqrt{y} / |D|) y^{(k+1)/2} (y - N) dy
\\&\qquad
  + \mathbbm{1}_{\sigma \le k - 5 / 2}
  \int_{N+1} ^{\infty} K_{k-2} (4 \pi \sqrt{y} / |D|)
      ((y - 1)/ y )^{\lfloor \sigma \rfloor}
\\&\qquad\qquad\qquad\qquad\qquad
      (y - 1)^{-k + 5 / 2} y^{3 k / 2 - 2} ( y - N) dy
\\&\qquad
  + \mathbbm{1}_{\sigma \ge k + 1/2}
  \int_{N+1} ^{\infty} K_{k-2} (4 \pi \sqrt{y} / |D|)
      (y / (y-1))^{\lceil \sigma \rceil}
\\&\qquad\qquad\qquad\qquad\qquad
      (y - 1)^{k+1} y^{(-k + 1) / 2} (y - N) dy
   \quad\Bigg)
,
\\&\qquad
  \widetilde{g}_D (n)
\\&:=
  2 (2 \pi / |D|)^{2 - k + 2 s} / (\Gamma(s) \Gamma(s - k + 2))
\\&\qquad
  \cdot
  \int_n ^\infty K_{k-2} (4 \pi \sqrt{y} / |D|)
      \left((y/n)^{3 k / 2 - 2 - \lfloor \sigma \rfloor} + (y/n)^{\lceil \sigma \rceil - k/2}
      \right)
.
\end{align*}
Note that this bound is far away from the optimal one. So we also use a second estimate
\begin{align*}
  \tilde \eta_f (P, N, P', N', D, s)
=
  |\tilde Z_{f,D} (s, P, N) - \tilde Z_{f,D} (s, P', N')|
\end{align*}
with suitable $P',\,N'$. Later we will have to introduce a correction factor $\kappa$, since $\tilde \eta$  might get too small if the series converges slowly.
Nevertheless, the bounds $\kappa \cdot \tilde \eta_f$ obtained in this way are in general much tighter.

In the course of these calculations we will use $P = 1500,\,N = 7999$ and $P' = 1000,\,N' = 3999$ and neglect these parameters in notation. We will also omit the parameter $s$, if $s = k - 1$.

Denote the Fourier expansion of $f \in [\Gamma_2,\,\det ^k ]$ by 
\begin{align*}
f(z) = \sum_t a_f (t) \exp(2 \pi i \cdot \tr(t z))
.
\end{align*}
Here $\tr(\cdot)$ is the usual trace function. We will use Andrianov's formula for the Hecke eigenvalues (cf. \cite{kohnenkuss_boechererconjecture} for an easy-to-read statement). To do so, we need to assume that $a_f \left(\begin{smallmatrix}1 & 1/2 \\ 1/2 & 1\end{smallmatrix}\right) \ne 0$. This assumption holds for all forms we consider.

We have to calculate the fourier expansion of forms $f$ up to discriminants $-\sqrt{8000}$. This is no problem at all. But moreover we have to calculate certain fourier coefficients of discriminant up to $-3 \cdot 1499^2 \approx -6.75 \cdot 10^6$.
Namely, we have to calculate the fourier coefficients associated to $\left(\begin{smallmatrix}1 & p/2 \\ p/2 & p^2\end{smallmatrix}\right)$. This actually is a very hard problem.

We briefly describe the technical details of out method. By the results of the preceding section we can use a basis which consists of products of at most two elements of the Maass spezialschar in all considered cases.
We have to compute $\eta^{-6}$. Notice that this is the most expensive part of the computation of the lifts' Fourier expansion and we precalculate it.
Multiplication of such big polynomials with only positive coefficients using Fourier transforms consumes too much memory.
Here, we combine a multimodular approach with ordinary $\mathcal{O} (n^2)$ multiplication, which admits easy hard drive caching.

Assume that we have calculated the elements of the Maass spezialschar. Let us now focus on a non trivial product $g = g_1 \cdot g_2$. We have to calculate certain Fourier coefficients $a_g (t) = \sum_{t_1 + t_2 = t} a_{g_1}(t_1) \cdot a_{g_2}(t_2)$. Since $a_{g_1}$ and $a_{g_2}$ depend only on $d_i := \det t_i,\, \epsilon_i := \gcd(t_i)$, we can write
\begin{align*}
  a_g (t)
&=
  \sum_{\epsilon_1, \epsilon_2, d_1, d_2}
  v(\epsilon_1, \epsilon_2, d_1, d_2)
  \, a_{g_1} (\epsilon_1, d_1)
  \, a_{g_2} (\epsilon_2, d_2)
.
\end{align*}
The function $v : \ZZ^4 \rightarrow \NN$ can be precalculated. The formula above helps avoiding multiplying the same numbers multiple times and hence saves a lot of time.
Consult the code, which is available on the author's homepage, for more details.

We calculate the central values of twisted spinor L-series associated to cuspidal Hecke eigenforms which do not belong to the Maass spezialschar up to even weight $k < 40$ and for all discriminants $0 > D > -300$. The relevant quotients are denoted by $c'_f (D) := \tilde Z_{f,D} (s, P, N) / B_f (D)$.
The factor $B_f (D)$ is essentially the right hand side of the Böcherer conjecture as it has been defined in the introduction.
The results confirm the Böcherer conjecture. However, for high discriminants they are not very precise. A brief summary is given in table \ref{boechererconjecturecfvalues} and \ref{boechererconjecturecfquotients}. For the complete data download the according files on the author's homepage.
A precise statement is:
\begin{theorem}
\label{boechererconjectureverificationtheorem}
Let $20 \le k < 40$ and let $f$ be a weight $k$ cuspidal Hecke eigenform which is not element of the Maass spezialschar. Then the intervals
\begin{align*}
[c'_f (D) - \eta_f(D) / B_f (D),\, c'_f (D) + \eta_f (D) / B_f (D)]
\end{align*}
for discriminants $0 > D > -300$ have non empty intersection. Choosing $\kappa_f$ according to table \ref{boechererconjectureverificationkappa} the intervals
\begin{align*}
[c'_f (D) - \kappa_f \tilde \eta_f (D) / B_f (D),\, c'_f (D) + \kappa_f \tilde \eta_f (D) / B_f (D)]
\end{align*}
have non empty intersection, too.
\end{theorem}

One might also try to verify Deligne's conjecture \cite{deglinelfunctions} to some extend. But apart from the central value the convergence at critical points is very bad and at the same time the values grow very rapidly.
For example the formula given above yields an error bound $\eta_f (-3, 21) \approx -8.7 \cdot 10^{169}$ and the L-function's value is about $\widetilde Z_{f, -3} (21) \approx 4.55 \cdot 10^{168} \approx \tilde \eta_f (-3, 21)$ for $k = 20,\, D = -3,\, s = 21$. To overcome this obstacle one might consider lower weight forms for congruence subgroups.


\section{Two further applications}
\label{sec:applications}

As mentioned in the introduction one motivation to find a basis as given in section \ref{sec:generatingsets} was to provide computational access to high weights.
Indeed, on an ordinary Sun server we are able to compute the Fourier expansions of the basis up to a precision we need to apply Hecke operators in few seconds.
In this section we present two applications which become amenable by the new bases.

\subsection{The minimal discriminant of pivot sets}
One often encounters modular forms in other fields of mathematics. In order to use modular forms as a tool of investigation of these fields it is crucial to express them in terms of a given basis. This is often done by comparing Fourier coefficients. Occasionally, it is very expensive to calculate the Fourier expansions. Since the costs of such a computation grow rapidly in terms of the discriminant of the Fourier index, it is natural to ask the following question. Up to which discriminant do one has to calculate the fourier expansion of a Siegel modular form of degree 2 to  uniquely determine it.

More precisely, we fix $k$ and denote the Fourier expansion of a form $f$ by $\sum a_f (t) \exp(2 \pi i \cdot \tr(t z))$. Now we consider the minimal discriminant $D(p) = \min_i (-4 \det t_i)$ associated to any set of pivot indices $p = \{t_i\}_i$ of $[\Gamma_2,\,\det ^k]$.
By pivot indices we mean a set of indices $p = \{t_i\}_i$ such that the rank of the matrix $(a_f (t_i))_{i, f} $ where $f$ runs over a basis of $[\Gamma_2,\, \det ^k]$ equals $\dim [\Gamma_2,\, \det ^k]$.
Now the question can be formulated as follows. Can we express $\max_p D(p)$ in terms of $k$?

\begin{figure}[H]
\setlength{\unitlength}{0.240900pt}
\ifx\plotpoint\undefined\newsavebox{\plotpoint}\fi
\sbox{\plotpoint}{\rule[-0.200pt]{0.400pt}{0.400pt}}%
\begin{picture}(1650,809)(0,0)
\sbox{\plotpoint}{\rule[-0.200pt]{0.400pt}{0.400pt}}%
\put(111,111){\makebox(0,0){0}}
\put(173.0,111.0){\rule[-0.200pt]{4.818pt}{0.400pt}}
\put(173.0,205.0){\rule[-0.200pt]{4.818pt}{0.400pt}}
\put(111,300){\makebox(0,0){2}}
\put(173.0,300.0){\rule[-0.200pt]{4.818pt}{0.400pt}}
\put(173.0,394.0){\rule[-0.200pt]{4.818pt}{0.400pt}}
\put(111,488){\makebox(0,0){4}}
\put(173.0,488.0){\rule[-0.200pt]{4.818pt}{0.400pt}}
\put(173.0,582.0){\rule[-0.200pt]{4.818pt}{0.400pt}}
\put(111,676){\makebox(0,0){6}}
\put(173.0,676.0){\rule[-0.200pt]{4.818pt}{0.400pt}}
\put(439,41){\makebox(0,0){2}}
\put(439.0,82.0){\rule[-0.200pt]{0.400pt}{4.818pt}}
\put(791.0,82.0){\rule[-0.200pt]{0.400pt}{4.818pt}}
\put(1143,41){\makebox(0,0){4}}
\put(1143.0,82.0){\rule[-0.200pt]{0.400pt}{4.818pt}}
\put(1494.0,82.0){\rule[-0.200pt]{0.400pt}{4.818pt}}
\put(193.0,102.0){\rule[-0.200pt]{0.400pt}{160.921pt}}
\put(193.0,102.0){\rule[-0.200pt]{338.946pt}{0.400pt}}
\put(1600.0,102.0){\rule[-0.200pt]{0.400pt}{160.921pt}}
\put(193.0,770.0){\rule[-0.200pt]{338.946pt}{0.400pt}}
\put(369,723){\makebox(0,0){$\log(\min D(p))$}}
\put(1406,143){\makebox(0,0){$\log (k)$}}
\sbox{\plotpoint}{\rule[-0.500pt]{1.000pt}{1.000pt}}%
\put(223,111){\circle*{12}}
\put(366,111){\circle*{12}}
\put(467,111){\circle*{12}}
\put(546,242){\circle*{12}}
\put(610,242){\circle*{12}}
\put(664,242){\circle*{12}}
\put(711,263){\circle*{12}}
\put(752,263){\circle*{12}}
\put(789,353){\circle*{12}}
\put(823,353){\circle*{12}}
\put(854,372){\circle*{12}}
\put(882,372){\circle*{12}}
\put(908,378){\circle*{12}}
\put(932,425){\circle*{12}}
\put(955,425){\circle*{12}}
\put(976,440){\circle*{12}}
\put(996,449){\circle*{12}}
\put(1015,449){\circle*{12}}
\put(1033,478){\circle*{12}}
\put(1050,478){\circle*{12}}
\put(1067,490){\circle*{12}}
\put(1082,498){\circle*{12}}
\put(1097,503){\circle*{12}}
\put(1112,519){\circle*{12}}
\put(1126,519){\circle*{12}}
\put(1139,529){\circle*{12}}
\put(1152,537){\circle*{12}}
\put(1164,541){\circle*{12}}
\put(1176,553){\circle*{12}}
\put(1187,553){\circle*{12}}
\put(1199,562){\circle*{12}}
\put(1209,569){\circle*{12}}
\put(1220,572){\circle*{12}}
\put(1230,582){\circle*{12}}
\put(1240,582){\circle*{12}}
\put(1250,589){\circle*{12}}
\put(1259,596){\circle*{12}}
\put(1268,598){\circle*{12}}
\put(1277,607){\circle*{12}}
\put(1286,607){\circle*{12}}
\put(1294,614){\circle*{12}}
\put(1303,619){\circle*{12}}
\put(1311,621){\circle*{12}}
\put(1319,629){\circle*{12}}
\put(1326,629){\circle*{12}}
\put(1334,635){\circle*{12}}
\put(1341,640){\circle*{12}}
\put(1348,642){\circle*{12}}
\put(1356,648){\circle*{12}}
\put(1363,648){\circle*{12}}
\put(1369,654){\circle*{12}}
\put(1376,659){\circle*{12}}
\put(1383,660){\circle*{12}}
\put(1389,666){\circle*{12}}
\put(1395,666){\circle*{12}}
\put(1402,672){\circle*{12}}
\put(1408,676){\circle*{12}}
\put(1414,677){\circle*{12}}
\put(1420,683){\circle*{12}}
\put(1426,683){\circle*{12}}
\put(1431,687){\circle*{12}}
\put(1437,692){\circle*{12}}
\put(1442,692){\circle*{12}}
\put(1448,698){\circle*{12}}
\put(1453,699){\circle*{12}}
\put(1459,702){\circle*{12}}
\put(1464,706){\circle*{12}}
\put(1469,707){\circle*{12}}
\put(1474,712){\circle*{12}}
\put(1479,713){\circle*{12}}
\put(1484,716){\circle*{12}}
\put(1489,719){\circle*{12}}
\put(1493,720){\circle*{12}}
\put(1498,724){\circle*{12}}
\put(1503,726){\circle*{12}}
\put(1507,728){\circle*{12}}
\put(1512,732){\circle*{12}}
\put(1516,732){\circle*{12}}
\put(1521,737){\circle*{12}}
\put(1525,738){\circle*{12}}
\put(1530,740){\circle*{12}}
\put(1534,744){\circle*{12}}
\put(1538,744){\circle*{12}}
\put(1542,748){\circle*{12}}
\put(1546,749){\circle*{12}}
\sbox{\plotpoint}{\rule[-0.200pt]{0.400pt}{0.400pt}}%
\put(193.0,102.0){\rule[-0.200pt]{0.400pt}{160.921pt}}
\put(193.0,102.0){\rule[-0.200pt]{338.946pt}{0.400pt}}
\put(1600.0,102.0){\rule[-0.200pt]{0.400pt}{160.921pt}}
\put(193.0,770.0){\rule[-0.200pt]{338.946pt}{0.400pt}}
\end{picture}
\caption{\label{minmaxdiscriminants_plot}Double logarithmic plot of maximal minimal discriminants for weights up to $172$.}
\end{figure}
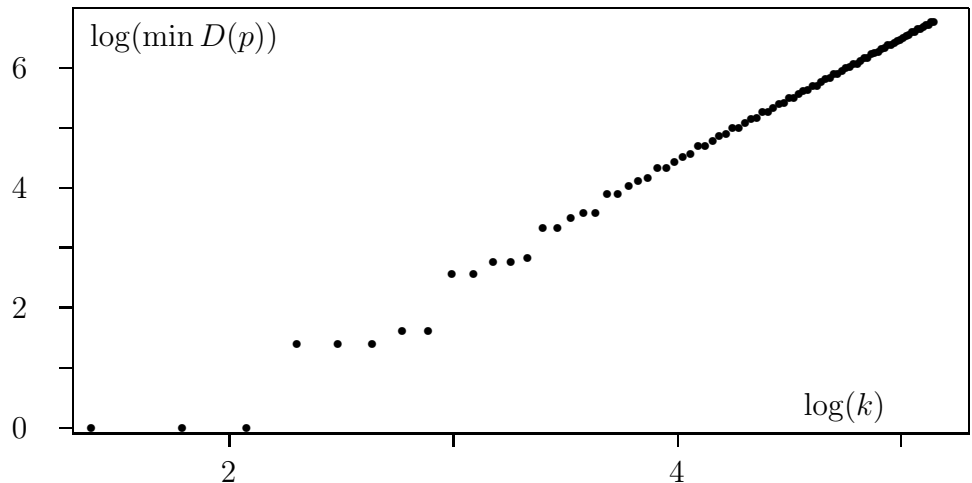

Using the basis presented in the preceding section, we can give a partial answer to this question. The double logarithmic plot \ref{minmaxdiscriminants_plot} of the data in table \ref{minmaxdiscriminants} reveals that for high weights $k$ the slope tends to $2$. This supports the following conjecture.
\begin{conjecture}
Fix $k$ and let $p$ run over all sets of pivot indices associated to weight $k$ forms. Then
\begin{align*}
  \max_p D(p)
=
  - c k^2 + \mathcal{O}(1) \quad \text{as }k\rightarrow \infty
,
\end{align*}
where $c \approx 736 / 158^2 < 0.03$.
\end{conjecture}

Note that due to the Victor-Miller construction the analogous result for $[\Gamma_1, \det ^k]$ is $\max_p D(p) = -k$.
In a recent preprint Poor and Yuen show a similar result implying $c \le 8 / 225 \approx 0.035$ (cf. \cite{pooryuenparamoduluar}).

\subsection{Irreducibility of the Hecke action}

Let $K : [\Gamma_1, k] \rightarrow [\Gamma_2, k]$ be the Klingen-Eisenstein lift.
The Hecke invariant splitting of the spaces of Siegel modular forms
\begin{align*}
[\Gamma_2, k] = K [\Gamma_1,k]_0 \oplus V [\Gamma_J, 1, k] \oplus S'_k 
\end{align*}
is well known.

In analogy to the degree $1$ case one might expect, that $S' _k$ is irreducible with respect to the rational Hecke action for all $k$.
But computations in \cite{skoruppasiegelcomputations} revealed that $S'_k$ is reducible if $k \in \{24,\, 26\}$. He found the aforementioned eigenforms $\Upsilon_{k a},\,\Upsilon_{k b}$.
Nevertheless, the conjecture that the Hecke action on spaces of Siegel modular forms should be ``as irreducible as possible'' is commonly supposed to hold true.
Using the computational benefit of the basis presented in section \ref{sec:generatingsets} we were able to prove
\begin{theorem}
The spaces $S'_k$ are irreducible with respect to the Hecke action for all even $k$ within the range $28 \le k \le 150$.
\end{theorem}
To prove this the author calculated the matrix associated to the Hecke operator $T(2)$ with respect to our basis using Sage and calculated its minimal polynomial using Magma.

\subsection{Interesting properties of the basis}
There are surprising properties of the bases presented in table \ref{spanningmaassproductstable}. These are canonical in the sense that the maximum of the weights of all lifts involved in non trivial products is minimal within the set of all possible bases of that kind.

We can read off the noteworthy fact that $14|_1 \cdot 16|_1$ is contained in the $6$ dimensional space spanned by the Maass spezialschar and the three elements $14|_0 \cdot 16|_0,\, 14|_0 \cdot 16|_1,\,14|_0 \cdot 16|^1$, whereas $\dim [\Gamma_2, \det ^{30}] = 11$. Hence the products we use to span the spaces of Siegel modular forms satisfy more linear relations then one might expect in general.

Next, we consider the two Hecke eigenforms $\Upsilon_{24a}$ and $\Upsilon_{24b}$ of weight $24$ normalized as in \cite{skoruppasiegelcomputations}. Expressing them in terms of the given basis yields two vectors $v_a,\, v_b \in \QQ^{8}$.
The fact that $(v_a)_i = 0 \Leftrightarrow (v_b)_i = 0$ as one can see in table \ref{upsiloncoordinatestable24} is surprising. An analogous result holds for the Hecke eigenforms $\Upsilon_{26a}$ and $\Upsilon_{26b}$ (see table \ref{upsiloncoordinatestable26}).

There is one further phenomenon which presumably is of arithmetic interest. The numerators as well as denominators of the entries in these vectors tend to be smooth.
Notice that $\Upsilon_{24 *}$ and $\Upsilon_{26 *}$ are canonically normalized in the sense that they are primitive integral eigenforms.

\vspace{1em}
Let us remark that the framework used to calculate this basis as well as the applications are available on the author's homepage. The framework is a branch of the tools the author and his coauthors present in \cite{computationofsiegelmodularforms} aiming at a greater flexibility of the code.


\section{Tables and plots}


\begin{table}[H]
\begin{tabular}{ll}
weight & products \\
\midrule[0.1em]
$20$   & $10|_0 \cdot 10|_0,\, 10|_0 \cdot 10|^1$  \\
$22$   & $10|_0 \cdot 12|_0,\, 10|_0 \cdot 12|_1$  \\
$24$   & $12|_0 \cdot 12|_0,\, 12|_0 \cdot 12|_1,\,%
          12|_1 \cdot 12|_1,\, 10|_0 \cdot 14|_0$   \\
$26$   & $12|_0 \cdot 14|_0,\, 12|_0 \cdot 14|^1,\,%
          12|_1 \cdot 10|_0$   \\
$28$   & $14|_0 \cdot 14|_0,\, 14|_0 \cdot 14|^1,\,%
          14|^1 \cdot 14|^1,\, 12|_0 \cdot 16|_0,\,%
          12|_0 \cdot 16|_1$   \\
$30$   & $14|_0 \cdot 16|_0,\, 14|_0 \cdot 16|_1,\,%
          14|_0 \cdot 16|^1,\, 14|_1 \cdot 16|_0,\,%
          14|^1 \cdot 16|_1,\, 12|_0 \cdot 18|_0$   \\
$32$   & $16|_0 \cdot 16|_0,\, 16|_0 \cdot 16|_1,\,%
          16|_0 \cdot 16|^1,\, 16|_1 \cdot 16|_1,\,%
          16|_1 \cdot 16|^1,\, 16|^1 \cdot 16|^1$   \\
       & $14|_0 \cdot 18|_0$ \\
$34$   & $16|_0 \cdot 18|_0,\, 16|_0 \cdot 18|_1,\,%
          16|_0 \cdot 18|^1,\, 16|_1 \cdot 18|_1,\,%
          16|_1 \cdot 18|^1,\, 16|_1 \cdot 18|^1$   \\
       & $16|^1 \cdot 18|_0,\, 14|_0 \cdot 20|_0$ \\
\end{tabular}
\caption{\label{spanningmaassproductstable}Products generating the spaces of Siegel modular forms, ignoring the Maass spezialschar}
\end{table}

\begin{table}[H]
\begin{tabular}{rr}
$\Upsilon_{24a}$   & $\Upsilon_{24b}$ \\
\midrule[0.1em]
$0$                & $0$              \\
$-29^2 \cdot 2237/2^5 \cdot 3^4 \cdot 5 \cdot 7 \cdot 13$ &
$416761/2^3 \cdot 3^5 \cdot 5^2 \cdot 7$  \\
$-52956193/2^2 \cdot 3^3 \cdot 5 \cdot 7 \cdot 13$  &
$-11 \cdot 83 \cdot 4987/3^4 \cdot 5^2 \cdot 7$ \\
$-227 \cdot 2969/2^5 \cdot 3^4 \cdot 5 \cdot 7 \cdot 13$  &
$937 \cdot 947/2^3 \cdot 3^5 \cdot 5^2 \cdot 7$  \\
$0$                & $0$              \\
$11 \cdot 157/2^6 \cdot 3^5 \cdot 5 \cdot 7 \cdot 13$     &
$13 \cdot 83/2^4 \cdot 3^6 \cdot 5^2 \cdot 7$   \\
$2^2 \cdot 5 \cdot 11 \cdot 157 \cdot 661/3$       &
$-2^4 \cdot 5 \cdot 7 \cdot 13^3 \cdot 83/3^2$   \\
$0$                & $0$              
\end{tabular}
\caption{\label{upsiloncoordinatestable24}Coordinates of the exceptional eigenforms of weight $24$ with respect to the basis given in table \ref{spanningmaassproductstable}.}

\end{table}

\begin{table}[H]
\begin{tabular}{rr}
$71 \cdot 139 \cdot \Upsilon_{26a}$          & $71 \cdot 139 \cdot \Upsilon_{26b}$    \\
\midrule[0.1em]
$0$                       & $0$                 \\
$10718579/2^3 \cdot 13^2$       &
$-3 \cdot 5 \cdot 54059/2^3$     \\
$61 \cdot 79 \cdot 3967087/2^4 \cdot 3 \cdot 5^2 \cdot 7 \cdot 13^2$ &
$-2251 \cdot 9923/2^4 \cdot 3 \cdot 5$ \\
$3 \cdot 1703358596089/2 \cdot 5^2 \cdot 7 \cdot 13^2$ &
$-3 \cdot 191 \cdot 10784339/2 \cdot 5$ \\
$0$                       & $0$                 \\
$-11 \cdot 29 \cdot 839/2 \cdot 5^2 \cdot 7 \cdot 13^2$       &
$4177/2^3 \cdot 3 \cdot 5$      \\
$11 \cdot 29 \cdot 37 \cdot 167/2^5 \cdot 3 \cdot 7 \cdot 13$        &
$-5 \cdot 7 \cdot 13 \cdot 19 \cdot 37/2^5 \cdot 3$
\end{tabular}
\caption{\label{upsiloncoordinatestable26}Coordinates of the exceptional eigenforms of weight $26$ with respect to the basis given in table \ref{spanningmaassproductstable}.}
\end{table}

\begin{table}[H]
\begin{tabular}{lrrrrrrrr}
weight     & $20$  & $22$ & $24a$ & $24b$ & $26a$ & $26b$ & $28$ & $30$  \\ 
\midrule[0.1em]
$\kappa_f$ & $180$ & $90$ & $70$  & $45$  & $400$ & $220$ & $90$ & $150$ \\
\midrule[0.2em]
weight     & $32$  & $34$   & $36$ & $38$ \\
\midrule[0.1em]
$\kappa_f$ & $120$ & $2300$ & $85$ & $400$
\end{tabular}
\caption{\label{boechererconjectureverificationkappa}Correction factors $\kappa$ in theorem \ref{boechererconjectureverificationtheorem}}
\end{table}

\begin{table}[H]
\begin{tabular}{lrlll}
weight & $D$ & $c'_f(D)$ & $\eta(D) / B_f (D) $ & $\tilde \eta (D) / B_f (D)$ \\
\midrule[0.1em]
$20$ & $-\hphantom{0}3$  & $2.067215202868765 \cdot 10^{11}$ 
      & $8.7123 \cdot 10^{-41}$           & $2.4462 \cdot 10^{-31}$ \\
     & $-\hphantom{0}4$  & $2.067215202868765 \cdot 10^{11}$
      & $2.6633 \cdot 10^{-24}$           & $2.9904 \cdot 10^{-17}$ \\
     & $-15$             & $2.067220228408808 \cdot 10^{11}$
      & $8.2915 \cdot 10^{\hphantom{-}25}$  & $1.8868 \cdot 10^{\hphantom{-0}6}$ \\
     & $-20$             & $2.066999728524787 \cdot 10^{11}$
      & $3.3599 \cdot 10^{\hphantom{-}33}$  & $2.2640 \cdot 10^{\hphantom{-0}8}$ \\
\midrule[0.05em]
$28$ & $-\hphantom{0}3$  & $5.105310601780946 \cdot 10^{14}$
      & $5.6429 \cdot 10^{-32}$          & $3.3222 \cdot 10^{-22}$ \\
     & $-\hphantom{0}4$  & $5.105310601780946 \cdot 10^{14}$
      & $9.7141 \cdot 10^{-17}$          & $3.4976 \cdot 10^{-10}$ \\
     & $-15$             & $5.105287106899537 \cdot 10^{14}$
      & $2.0716 \cdot 10^{\hphantom{-}45}$ & $1.2487 \cdot 10^{\hphantom{-}11}$ \\
     & $-20$             & $5.108464527256237 \cdot 10^{14}$
      & $8.7051 \cdot 10^{\hphantom{-}53}$ & $2.0228 \cdot 10^{\hphantom{-}12}$ \\
\midrule[0.05em]
$30$ & $-\hphantom{0}3$  & $7.494659605198417 \cdot 10^{15}$
      & $6.3380 \cdot 10^{-30}$          & $3.3652 \cdot 10^{-21}$ \\
     & $-\hphantom{0}4$  & $7.494659605198417 \cdot 10^{15}$
      & $1.4559 \cdot 10^{-13}$          & $3.3121 \cdot 10^{-\hphantom{0}7}$ \\
     & $-15$             & $7.493428754840905 \cdot 10^{15}$
      & $4.6698 \cdot 10^{\hphantom{-}50}$ & $3.8657 \cdot 10^{\hphantom{-}12}$ \\
     & $-20$             & $7.525312838702049 \cdot 10^{15}$
      & $2.9997 \cdot 10^{\hphantom{-}60}$ & $2.7895 \cdot 10^{\hphantom{-}14}$ \\
\midrule[0.05em]
$32$ & $-\hphantom{0}3$  & $1.205163535250015 \cdot 10^{17}$
      & $6.3368 \cdot 10^{-28}$          & $3.5683 \cdot 10^{-18}$ \\
     & $-\hphantom{0}4$  & $1.205163535250015 \cdot 10^{17}$
      & $1.5997 \cdot 10^{-12}$          & $6.5727 \cdot 10^{-\hphantom{0}7}$ \\
     & $-15$             & $1.205732652809925 \cdot 10^{17}$
      & $3.2290 \cdot 10^{\hphantom{-}57}$ & $3.7921 \cdot 10^{\hphantom{-}15}$ \\
     & $-20$             & $1.205060112394051 \cdot 10^{17}$
      & $2.9225 \cdot 10^{\hphantom{-}64}$ & $7.7758 \cdot 10^{\hphantom{-}13}$ \\
\midrule[0.05em]
$34$ & $-\hphantom{0}3$  & $1.086138503038145 \cdot 10^{18}$
      & $5.6918 \cdot 10^{-26}$          & $2.5842 \cdot 10^{-17}$ \\
     & $-\hphantom{0}4$  & $1.086138503038145 \cdot 10^{18}$
      & $9.1579 \cdot 10^{-10}$          & $6.7761 \cdot 10^{-\hphantom{0}4}$ \\
     & $-15$             & $1.086049155614554 \cdot 10^{18}$
      & $1.3262 \cdot 10^{\hphantom{-}61}$ & $9.5845 \cdot 10^{\hphantom{-}15}$ \\
     & $-20$             & $1.084637644518613 \cdot 10^{18}$
      & $9.1134 \cdot 10^{\hphantom{-}70}$ & $3.1879 \cdot 10^{\hphantom{-}16}$ \\
\midrule[0.05em]
$36$ & $-\hphantom{0}3$  & $1.684497849534415 \cdot 10^{19}$
      & $4.6300 \cdot 10^{-24}$          & $5.5959 \cdot 10^{-16}$ \\
     & $-\hphantom{0}4$  & $1.684497849534415 \cdot 10^{19}$
      & $6.5481 \cdot 10^{-\hphantom{0}9}$ & $4.5361 \cdot 10^{-\hphantom{0}3}$  \\
     & $-15$             & $1.704113931864776 \cdot 10^{19}$
      & $4.0039 \cdot 10^{\hphantom{-}68}$ & $1.2502 \cdot 10^{\hphantom{-}18}$\\
     & $-20$             & $1.683266203530922 \cdot 10^{19}$
      & $3.3778 \cdot 10^{\hphantom{-}75}$ & $5.0592 \cdot 10^{\hphantom{-}16}$ \\
\midrule[0.05em]
$38$ & $-\hphantom{0}3$  & $6.357076618815059 \cdot 10^{19}$
      & $3.4351 \cdot 10^{-22}$          & $1.1738 \cdot 10^{-14}$ \\
     & $-\hphantom{0}4$  & $6.357076618815059 \cdot 10^{19}$
      & $1.8281 \cdot 10^{-\hphantom{0}7}$ & $4.3949 \cdot 10^{-\hphantom{0}1}$ \\
     & $-15$             & $6.517241676742821 \cdot 10^{19}$
      & $8.0534 \cdot 10^{\hphantom{-}73}$ & $2.9605 \cdot 10^{\hphantom{-}18}$ \\
     & $-20$             & $6.375510867491227 \cdot 10^{19}$
      & $1.6586 \cdot 10^{\hphantom{-}81}$ & $7.7408 \cdot 10^{\hphantom{-}17}$
\end{tabular}
\caption{\label{boechererconjecturecfvalues}Numerical estimates $c'_f (D)$ of the constant $c_f$ occurring in the Böcherer conjecture.}
\end{table}

\begin{table}[H]
\begin{tabular}{lrrrrrrr}
$D$ \textbackslash weight \ 
          & $20$     & $28$     & $30$     & $32$     & $34$     & $36$     
          & $38$       \\
\midrule[0.1em]
$-4$      & $0.0000$ & $0.0000$ & $0.0000$ & $0.0000$ & $0.0000$ & $0.0000$ 
          & $0.0000$ \\
$-7$      & $0.0000$ & $0.0000$ & $0.0000$ & $0.0000$ & $0.0000$ & $0.0000$ 
          & $0.0000$ \\
$-20$     & $0.0001$ & $0.0006$ & $0.0041$ & $0.0001$ & $0.0014$ & $0.0007$ 
          & $0.0029$ \\
$-103$    & $1.3565$ & $0.1592$ & $3.3038$ & $0.4526$ & $1.2287$ & $1.0546$ 
          & $1.9854$ \\
$-203$    & $0.2641$ & $0.4993$ & $2.2404$ & $0.1465$ & $0.6835$ & $0.6584$ 
          & $0.4231$ \\
$-299$    & $0.2152$ & $0.2511$ & $0.1197$ & $0.7238$ & $1.1443$ & $0.2187$ 
          & $0.2517$ \\
\end{tabular}
\caption{\label{boechererconjecturecfquotients}Some values of $|\log (c'_f (D) / c'_f (-3))|$.}
\end{table}

\begin{table}[H]
\begin{tabular}{rlrlrl}
weight & $\max D(p)$ & weight & $\max D(p)$ & weight & $\max D(p)$ \\
\midrule[0.1em]
 $100$ & $301$ &  $120$ & $433$ &  $140$ & $589$ \\
 $102$ & $301$ &  $122$ & $436$ &  $142$ & $596$ \\
 $104$ & $320$ &  $124$ & $456$ &  $144$ & $616$ \\
 $106$ & $337$ &  $126$ & $477$ &  $146$ & $641$ \\
 $108$ & $341$ &  $128$ & $481$ &  $148$ & $645$ \\
 $110$ & $364$ &  $130$ & $508$ &  $150$ & $676$ \\
 $112$ & $365$ &  $132$ & $513$ &  $152$ & $685$ \\
 $114$ & $385$ &  $134$ & $533$ &  $154$ & $705$ \\
 $116$ & $404$ &  $136$ & $556$ &  $156$ & $732$ \\
 $118$ & $408$ &  $138$ & $560$ &  $158$ & $736$
\end{tabular}
\caption{\label{minmaxdiscriminants}Maximal minimal discriminants of pivot sets for weights \mbox{$100 \le k < 160$}.}
\end{table}

\bibliographystyle{apalike}
\bibliography{bibliography}

\end{document}